\documentclass[12pt,a4paper]{amsart}

\usepackage {amsfonts}
\usepackage[english]{babel}
\def\g{\gamma}
\def\G{\Gamma}
\def\d{\delta}
\def\p{\varphi}
\def\e{\varepsilon}
\def\l{\lambda}
\def\L{\Lambda}
\def\la{\langle}
\def\ra{\rangle}
\def\R{{\mathbb R}}
\def\C{{\mathbb C}}
\def\N{{\mathbb N}}
\def\hT{\frak T}

\def\r{\rho}
\def\bs{~\hfill\rule{7pt}{7pt}}

\DeclareMathOperator{\supp}{supp}
\DeclareMathOperator{\diam}{diam}

\newtheorem{Th}{Theorem}
\newtheorem{Lem}{Lemma}

\begin{document}

\title{Uniqueness Theorems for Fourier Quasicrystals and Temperate Distributions with Discrete Support}

\author{S.Yu.Favorov}

\address{Sergii Favorov,
\newline\hphantom{iii}  Karazin's Kharkiv National University
\newline\hphantom{iii} Svobody sq., 4,
\newline\hphantom{iii} 61022, Kharkiv, Ukraine}
\email{sfavorov@gmail.com}

\maketitle {\small
\begin{quote}
\noindent{\bf Abstract.}
It is proved that if some points of the supports of two Fourier quasicrystals approach each other while tending to infinity and the same is true for the masses at these points, then these quasicrystals coincide. A similar statement is obtained for a certain class of discrete temperate distributions.
\medskip


\end{quote}
}

\section{introduction}
P.Kurasov and R.Suhr \cite{KS} noted that if zeros of two holomorphic almost periodic functions in a strip get closer at infinity, then the zero sets of these functions coincide. This result can be interpreted as the coincidence of two almost periodic discrete sets if they get closer at infinity. It is natural to expect the same effect for other almost periodic objects, in particular, for
Fourier quasicrystals or, in general, for distributions with discrete support and spectrum.
\bigskip

Denote by $S(\R^d)$ the Schwartz space of test functions $\p\in C^\infty(\R^d)$ with finite norms
\begin{equation}\label{N}
  N_m(\p)=\sup_{\R^d}(\max\{1,|x|\})^m\max_{\|k\|\le m} |(D^k\p)(x))|,\quad m=0,1,2,\dots,
\end{equation}
  $k=(k_1,\dots,k_d)\in(\N\cup\{0\})^d,\  \|k\|=\|k\|_\infty=\max\{k_1,\dots,k_d\},\ D^k=\partial^{k_1}_{x_1}\dots\partial^{k_d}_{x_d}$. These norms generate a topology on $S(\R^d)$, and elements of the space $S^*(\R^d)$
  of continuous linear functionals on $S(\R^d)$ are called temperate distributions.

The Fourier transform of a temperate distribution $f$ is defined by the equality
\begin{equation}\label{B}
\hat f(\p)=f(\hat\p)\quad\mbox{for all}\quad\p\in S(\R^d),
\end{equation}
where
$$
   \hat\p(y)=\int_{\R^d}\p(x)\exp\{-2\pi i\la x,y\ra\}dx
 $$
is the Fourier transform of the function $\p$.
\medskip

We will say that a distribution (or a measure) $f$ on $\R^d$  is {\it discrete}, if for each $\l\in\supp f$ there is $\e=\e(\l)>0$ such that $B(\l,\e)\cap\supp f=\{\l\}$,
and  {\it uniformly discrete}, if there is $\e>0$ such that $B(\l,\e)\cap B(\l',\e)=\emptyset$ for all $\l,\l'\in\supp f,\,\l\neq\l'$;
a measure $\mu$ is {\it atomic} if $\mu=\sum_{\l\in\L}a_\l\d_\l$ with $a_\l\in\C$ and countable $\L$, in this case we will write $a_\l=\mu(\l)$.

Here $B(x,r)=\{y\in\R^d:\,|y-x|<r\}$, and $\d_\l$ means the unit mass at the point $\l\in\R^d$.

A complex measure $\mu\in S^*(\R^d)$ is a {\it Fourier quasicrystal} if $\mu$ and its Fourier transform $\hat\mu$ are discrete  measures, and the measures $|\mu|$ and $|\hat\mu|$ belong to $S^*(\R^d)$.

Note that the condition $\mu\in S^*(\R^d)$ do not imply $|\mu|\in S^*(\R^d)$ (see \cite{LA} or Remark 3 of the present article).

Such measures are the main object in the theory of  Fourier quasicrystals (see \cite{L}, \cite{LO1}-\cite{Mo}, \cite{F2}). The corresponding notion was inspired by experimental discovery
 of non-periodic atomic structures with diffraction patterns consisting of spots, which was made in the mid '80s.

We will say that a complex measure $\mu$ is a {\it sparse} Fourier quasicrystal, when $\mu$ is discrete, $\mu\in S^*(\R^d)$, $\hat\mu$ is atomic, $|\hat\mu|\in S^*(\R^d)$,
 and numbers of elements $\#\{\supp\mu\cap B(x,1)\}$ are uniformly bounded in $x\in\R^d$.

Note that, compared with the classical definition of Fourier quasicrystal, we have weakened the conditions on the measure $\hat\mu$ and removed the requirement $|\mu|\in S^*(\R^d)$.

Clearly, a Fourier quasicrystal with a uniformly discrete support is a sparse Fourier quasicrystal.

\begin{Th}\label{T1}
If two sparse Fourier quasicrystals $\mu=\sum_{\l\in\L} \mu(\l)\d_{\l},\ \nu=\sum_{\g\in\G} \nu(\g)\d_{\g}$ under appropriate numbering $\L=\{\l_n\}_{n=1}^\infty,\,\G=\{\g_n\}_{n=1}^\infty$ have
the properties
\begin{equation}\label{a}
\l_n-\g_n\to0  \quad\mbox{ and }\quad  \mu(\l_n)-\nu(\g_n)\to0 \quad \mbox{ as }\ n\to\infty,
\end{equation}
 then the measures $\mu,\,\nu$  coincide.
\end{Th}

 The conditions of the theorem can be significantly weakened. First, the sparseness of measures and conditions (\ref{a}) can only be checked on a set $E$ with the property
 \begin{equation}\label{E}
 \exists\ \{B(x_k,r_k)\}_{k=1}^\infty,\quad r_k\to\infty,\quad r_k/|x_k|\to0\quad\mbox{such that}\quad\bigcup_{k=1}^\infty B(x_k,r_k)\subset E.
 \end{equation}
\begin{Th}\label{T1a}
Let $\mu=\sum_{\l\in\L} \mu(\l)\d_{\l},\ \nu=\sum_{\g\in\G} \nu(\g)\d_{\g}$ be discrete measures from $S^*(\R^d)$ such that
\begin{equation}\label{b3}
  \hat\mu, \ \hat\nu\quad\mbox{ are atomic measures, and}\quad |\hat\mu|,\ |\hat\nu|\in S^*(\R^d),
\end{equation}
 which  "sparse" on a set $E$ with  property (\ref{E}), i.e., there is $N<\infty$ such that for all $x\in E$
 $$
  \#\{n:\,\l_n\in B(x,1)\}\le N,\quad \quad\#\{n:\,\g_n\in B(x,1)\}\le N.
$$
If under appropriate numbering
\begin{equation}\label{a1}
\l_n-\g_n\to0  \quad\mbox{ and }\quad  \mu(\l_n)-\nu(\g_n)\to0
\end{equation}
as $n\to\infty$ and $\l_n$ or $\g_n$ belong to $E$, then the measures $\mu,\,\nu$ coincide.
\end{Th}

Second, points of the support of measures can be combined into groups that are located in a certain sense sparse and behave like points in (\ref{a1}). Moreover, the measures $\mu,\ \nu$ may be atomic:
\begin{Th}\label{T2}
Let $\mu,\,\nu$ be atomic measures from $S^*(\R^d)$,  whose Fourier transforms $\hat\mu, \ \hat\nu$ satisfy (\ref{b3}). Suppose that  there exist a set $E$ with  property (\ref{E}), disjoint
Borel sets $\L_n\subset E$, and  disjoint Borel sets $\G_n\subset E$ such that the restrictions of measures to $E,\ \L_n,\ \G_n$ satisfy the conditions
 $$
  \mu\mid_E=\sum_n\mu\mid_{\L_n},\qquad \nu\mid_E=\sum_n\nu\mid_{\G_n},
 $$
(we do not require $\L_n\neq\emptyset$ or $\G_n\neq\emptyset$). Also, for some $N<\infty$
 \begin{equation}\label{b1}
 \quad \#\{n:\,(\L_n\cup\G_n)\cap B(x,1)\neq\emptyset\}\le N,\quad \forall\ x\in E.
\end{equation}
If
\begin{equation}\label{b2}
\diam\{\L_n\cup\G_n\}\to0 \quad\mbox{ and }\quad  \mu(\L_n)-\nu(\G_n)\to0\quad\mbox{as}\quad n\to\infty,
\end{equation}
 then the measures $\mu,\,\nu$ coincide.
\end{Th}

{\bf Remark 1}. In the case of discrete measures $\mu,\ \nu$ and of one-point sets $\L_n,\,\G_n$ Theorem \ref{T2} coincides with Theorem \ref{T1a}.

\medskip

{\bf Remark 2}. Conditions $\mu,\,\nu\in S^*(\R^d)$  and (\ref{b3}) can be replaced by the following:  the functions
$$
\int\p(x-t)\mu(dx)\quad\mbox{ and }\quad \int\p(x-t)\nu(dx)
$$
are almost periodic for every $\p\in C^\infty$ with compact support.

\medskip
Under some additional conditions, the uniqueness theorem also holds for distributions with discrete supports.
 Note that by \cite{F3}, Proposition 3.1, for every temperate distribution $F$ with discrete support there is $m<\infty$ such that
$$
F=\sum_{\l\in\L}\sum_{\|j\|\le m}p_{\l,j}D^j\d_\l,\quad j\in(\N\cup\{0\})^d,\quad p_{\l,j}\in\C.
$$

\begin{Th}\label{T3}
Let
$$
f=\sum_{\l\in\L}\sum_{\|j\|\le m}p_{\l,j}D^j\d_\l,\quad g=\sum_{\g\in\G}\sum_{\|j\|\le m}q_{\g,j}D^j\d_{\g}
$$
be temperate distributions with discrete supports $\L,\,\G$ such that
\begin{equation}\label{f1}
  \hat f, \ \hat g\quad\mbox{ are atomic measures, and}\quad |\hat f|,\ |\hat g|\in S^*(\R^d),
\end{equation}
and there be a set $E$ with  property (\ref{E}), disjoint sets $\L_n$, and  disjoint sets $\G_n$  with properties (\ref{b1}) such that
$$
  \L\cap E=\bigcup_n\L_n,\qquad \G\cap E=\bigcup_n\G_n,
$$
(we do not require $\L_n\neq\emptyset$ or $\G_n\neq\emptyset$). If
\begin{equation}\label{f3}
\diam\{\L_n\cup\G_n\}\to0, \quad\quad  \sum_{\l\in\L_n} p_{\l,j}-\sum_{\g\in\G_n}q_{\g,j}\to0 \quad \forall j\quad\mbox{as}\quad n\to\infty,
\end{equation}
and
\begin{equation}\label{f4}
\max_j\sup_n\sum_{\l\in\L_n}|p_{\l,j}|<\infty,\quad \max_j\sup_n\sum_{\g\in\G_n}|q_{\g,j}|<\infty,
\end{equation}
then the distributions $f,\,g$ coincide.
\end{Th}

Measures and their Fourier transforms can be interchanged in Theorems \ref{T2} and \ref{T3}. Recall that the support of the Fourier transform of a measure $\mu$ is called a {\it spectrum} of $\mu$.

\begin{Th}\label{T4}
Let $\mu,\,\nu$ be atomic measures with discrete spectra $\tilde\L,\,\tilde\G$ respectively, $|\mu|,\,|\nu|\in S^*(\R^d)$,  and there be a set $E$ with  property (\ref{E}), disjoint sets $\tilde\L_n$,
and disjoint sets $\tilde\G_n$ such that
$$
\tilde\L\cap E=\bigcup_n\tilde\L_n,\qquad \tilde\G\cap E=\bigcup_n\tilde\G_n,\qquad \diam\{\tilde\L_n\cup\tilde\G_n\}\to0\quad\mbox{as}\quad n\to\infty,
$$
and conditions (\ref{b1}) are met.

 If either $\hat\mu,\,\hat\nu$ are measures, and
 $$
 \hat\mu(\tilde\L_n)-\hat\nu(\tilde\G_n)\to0\quad\mbox{as}\quad n\to\infty,
 $$
 or
 $$
\hat\mu=\sum_{\l\in\tilde\L}\sum_{\|j\|\le m}\tilde p_{\l,j}D^j\d_\l,\quad \hat\nu=\sum_{\g\in\tilde\G}\sum_{\|j\|\le m}\tilde q_{\g,j}D^j\d_{\g},
$$
and
$$
\sum_{\l\in\tilde\L_n}\tilde p_{\l,j}-\sum_{\g\in\tilde\G_n}\tilde q_{\g,j}\to0 \quad \forall j \quad\mbox{ as}\quad n\to\infty,
$$
$$
\max_j\sup_n\sum_{\l\in\tilde\L_n}|\tilde p_{\l,j}|<\infty,\quad \max_j\sup_n\sum_{\g\in\tilde\G_n}|\tilde q_{\g,j}|<\infty,
$$
  then the measures $\mu,\,\nu$ coincide.
\end{Th}

We also give an analogue of Theorem \ref{T1} for the Fourier quasicrystals:

\begin{Th}\label{T5}
If two Fourier quasicrystals $\mu,\,\nu$ with discrete "sparse" spectra $\tilde\L,\,\tilde\G$ (this means that
$$
\sup_{y\in\R^d}\#\{\tilde\L\cap B(y,1)\}<\infty,\quad \sup_{y\in\R^d}\#\{\tilde\G\cap B(y,1)\}<\infty)
$$
 under appropriate numbering $\tilde\L=\{\l_n\}_{n=1}^\infty,\,\tilde\G=\{\g_n\}_{n=1}^\infty$ have
the properties
$$
\l_n-\g_n\to0  \quad\mbox{ and }\quad  \hat\mu(\l_n)-\hat\nu(\g_n)\to0 \quad \mbox{ as }\ n\to\infty,
$$
 then the measures $\mu,\,\nu$  coincide.
\end{Th}

\section{auxiliary results}

\begin{Lem}\label{L1}
Let a positive measure $\mu$ belong to $S^*(\R^d)$. Then there is $N<\infty$ such that
\begin{equation}\label{R}
\mu(B(0,R))=O(R^N),
\end{equation}
 and for any Borel
function $H(x)$  such that $\sup_{x\in\R^d}|H(x)| (1+|x|^T)<\infty$  for all $T<\infty$ we get
$$
  \int_{\R^d}|H(x)|\mu(dx)<\infty.
$$
\end{Lem}

{\bf Proof of the Lemma}. Assume the converse. Then there is a sequence $R_n\to\infty$ such that $\mu(B(0,R_n))>R_n^n$. We may suppose that $R_{n+1}>2R_n$ for all $n$.
Take $\p(t)\in C^\infty(\R),\ 0\le\p(t)\le1$ such that $\p(t)=1$ for $t\le1$ and $\p(t)=0$ for $t\ge2$. Set
$$
\Psi(x)=\sum_n R_n^{-n}\p(|x|/R_n).
$$
Clearly, $\Psi\in C^\infty$ and
\begin{equation}\label{c}
  \int_{\R^d}\Psi(x)\mu(dx)\ge \sum_n R_n^{-n}\mu(B(0,R_n))=\infty.
\end{equation}
On there other hand, take any $K<\infty$ and $x$ such that $2R_{p-1}<|x|\le 2R_p$ with $p>K$. We have
$$
  |x|^K\Psi(x)=\sum_n |x|^K R_n^{-n}\p(|x|/R_n)<2^K R_p^{K-p}\sum_{n\ge p}R_p^p/R_n^n.
$$
Taking into account that $R_n>2^{n-p}R_p$ and $p\to\infty$ as $|x|\to\infty$, we obtain
$$
  |x|^K\Psi(x)<2^{K+1}R_p^{K-p}\to0 \quad \mbox{ as } |x|\to\infty.
$$
Similarly, one can check that $|x|^K\Psi^{(k)}(x)\to0$ for all $K$ and $k\in(\N\cup\{0\})^d$, therefore, $\Psi\in S(\R^d)$. Since $\mu\in S^*(\R^d)$, we get the contradiction with (\ref{c}).
Hence there exists $N$ such that $M(R):=\mu(B(0,R))\le C\max(1,R^N)$.

Furthermore, let $|H(x)|\le C_1|x|^{-N-1}$ for $|x|\ge1$.  Passing to polar coordinates and integrating in parts, we obtain
$$
  \int_{\R^d}|H(x)|\mu(dx)\le C_0+C_1\int_{|x|>1}|x|^{-N-1}\mu(dx)=C_0+C_2\int_1^\infty r^{-N-1}M(dr)
$$
$$
=C_0+C_2\left(\lim_{R\to\infty}\frac{M(R)}{R^{N+1}} -M(1)+(N+1)\int_1^\infty \frac{M(r)}{r^{N+2}}dr\right)<\infty.
$$
Lemma is proved. \bs

{\bf Remark 3}. Let
$$
 \mu=\sum_{n\in\N}2^{n-1}[\d_{n+2^{-n}}-\d_{n-2^{-n}}]
$$
be a measure on $\R$. For any $\p\in S(\R)$ we have for some $t_n\in[n-2^{-n},n+2^{-n}]$
$$
   \int_\R \p(t)\mu(dt)=\sum_{n\in\N}2^{n-1}[\p(n+2^{-n})-\p(n-2^{-n})]=\sum_{n\in\N}\p'(t_n),
$$
and
$$
  |(\mu,\p)|=\left| \int_\R \p(t)\mu(dt)\right|\le\sum_{n\in\N}t_n^{-2}N_2(\p)\le CN_2(\p),
$$
where $N_2(\p)$ is defined in (\ref{N}). Therefore, $\mu\in S^*(\R)$. On the other hand,
$$
  |\mu|(-n-2^{-n},n+2^{-n})=\sum_{1\le j\le n}2^{j+1}=2^{n+2},
$$
hence by Lemma \ref{L1}, the condition $|\mu|\in S^*(\R)$ contradicts to (\ref{R}).
\medskip

The proofs of our theorems are also based on the properties of almost periodic functions and distributions. Recall some definitions
related to the notion of almost periodicity \newline (a detailed exposition of the theory of almost periodic functions on $\R$ see, for example, in \cite{B} and \cite{LZ}, most of the results can
easily be generalized to functions on $\R^d$; almost periodic measures and distributions were introduced in \cite{LA} and \cite{R}, see also \cite{M1}, \cite{M2}, \cite{FK1}, \cite {F3}).

A set $A\subset\R^d$ is relatively dense, if there is $R<\infty$ such that every ball of radius $R$ intersects with $A$.

A continuous function $f$ on $\R^d$ is {\it almost periodic}, if
for every  $\e>0$ the set of $\e$-almost periods of $f$
  $$
  \{\tau\in\R^d:\,\sup_{t\in\R^d}|f(t+\tau)-f(t)|<\e\}
  $$
  is relatively dense in $\R^d$.

For example, for arbitrary $s_n\in\R^d$ the function
$$
    f(t)=\sum_n a_n e^{2\pi i\la t,s_n\ra}
$$
is almost periodic under the condition $\sum_n |a_n|<\infty$.

It was proved in \cite{B} that a finite family $\{f_j\}_{j=1}^M$ of almost periodic functions on $\R$ has a common relatively dense set
of  $\e$-almost periods for every $\e$. The same result for almost periodic functions on Euclidean spaces follows immediately from Bochner's criterion: a function $f(x),\,x\in\R$ is almost periodic if and only if
for every sequence $x_n$ there is a subsequence $x_{n'}$ such that the functions $f(x+x_{n'})$ converges uniformly in $x\in\R$. Its proof in \cite{LZ} practically without changes is
transferred to functions on $\R^d$ and even to mappings from $\R^d$ to $\R^M$.

 Now, if each $f_j$ satisfies  Bochner's criterion,  then the mapping $F=(f_1(x),\dots,f_M(x))$ satisfies this criterion too.
 It remains to notice that every $\e$-almost period of $F$ is an $\e$-almost period of every $f_j$.

Next, put $\p_t=\p(x-t)$ for any function $\p$ on $\R^d$.

A  measure $\mu$ is {\it almost periodic}, if the function $F(t)=\int\psi_t(x)\mu(dx)$
is almost periodic for any continuous function $\psi(x)$ with compact support.

A temperate distribution  $f$ is {\it almost periodic}, if for every $\p\in S(\R^d)$ the function $F(t)=f(\p_t)$
is almost periodic.

It is easy to prove that every almost periodic measure is an almost periodic temperate distribution.  Also, if $\mu$ is a positive measure, or satisfies the condition
$\sup_{x\in\R^d}|\mu(B(x,1))|<\infty$, then almost periodicity of $\mu$ in the sense of distributions implies almost periodicity in the sense of measures.
But there are discrete measures that almost periodic temperate distributions and
 not almost periodic measures   (\cite{LA}, \cite{M1}, \cite{FK2}).

\begin{Lem}\label{L2}
If $f\in S^*(\R^d)$, $\hat f$ is an atomic measure, and $|\hat f|\in S^*(\R^d)$, then $f$ is an almost periodic distribution.

 In particular, every Fourier quasicrystal is an almost periodic distribution.
\end{Lem}

{\bf Proof of the Lemma}.
Let $\hat f=\sum_n b_n\d_{s_n}$, then $|\hat f|=\sum_n |b_n|\d_{s_n}$. By (\ref{B}),  we have for each $\p\in S(\R^d)$
\begin{equation}\label{d}
  f(\p_t)=\int\check\p(y)e^{2\pi i\la t,y\ra}\hat  f(dy)=\sum_n b_n\check\p(s_n)e^{2\pi i\la t,s_n\ra},
\end{equation}
where $\check\p(y)e^{2\pi i\la t,y\ra}$ is the inverse Fourier transform of the function $\p(x-t)$.  Since $\check\p(y)\in S(\R^d)$, we  apply Lemma \ref{L1} and get
$$
  \sum_n |b_n||\check\p(s_n)|=\int|\check\p(y)||\hat f|(dy)<\infty.
$$
 Therefore the function $f(\p_t)$ is almost periodic.  \bs
\medskip

\section{proofs of the theorems}

{\bf Proof of Theorem \ref{T2}}. Assume the contrary $\mu\not\equiv\nu$. Since $\mu,\,\nu$ are atomic measures, we can find  $a\in\R^d$ such that $\nu(a)\neq\mu(a)$.
 Set
$$
\e=\frac{|\mu(a)-\nu(a)|}{3}.
$$
 Since
$$
 (|\mu|+|\nu|)(\{x:\,0<|x-a|<\r\})\to0\quad\mbox{as}\quad \r\to0,
  $$
 we see that for some $\r\le 1/2$
\begin{equation}\label{b}
 |\mu|(\{x:\,0<|x-a|<\r\})+|\nu|(\{x:\,0<|x-a|<\r\})<\e.
\end{equation}
Take $\p\in C^\infty(\R^d)$ such that $0\le\p\le1$,  $\supp\p\subset B(2)$, and $\p(x)=1$ for $|x|<1$.
 Let $N$ be a number from (\ref{b1}). Set for $j=1,\dots,2N+1$
$$
f_j(t)=\int\p\left(\frac{x-t}{2^{-j}\r}\right)\mu(dx),\quad  g_j(t)=\int\p\left(\frac{x-t}{2^{-j}\r}\right)\nu(dx), \quad H_j(t)=f_j(t)- g_j(t).
$$
Using (\ref{b3}) and applying Lemma \ref{L2}, we obtain that all the functions $f_j(t)$ and $g_j(t)$
 are almost periodic, and the functions  $H_j(t)$ too. Also, it follows from (\ref{b}) that
$$
|H_j(a)|>2\e\quad\forall\ j=1,\dots,2N+1.
$$
Denote by $\hT$ the set of all common $\e$-almost periods of the functions $H_j(t)$. We get
\begin{equation}\label{e}
|H_j(a+\tau)|>\e \quad \forall \tau\in\hT,\quad j=1,\dots,2N+1.
\end{equation}
Since $\hT$ is relatively dense, it follows from (\ref{E}) that for every $k>k_1$ there is $\tau_k\in\hT$ such that $B(x_k,r_k)\supset B(a+\tau_k,\r)$.
By (\ref{b1}), the set
$$
\{n:\,(\L_n\cup\G_n)\subset B(x_k,r_k)\}
$$
is discrete, hence, we have
\begin{equation}\label{n}
\min\{n:\,(\L_n\cup\G_n)\subset B(x_k,r_k)\}\to\infty\quad\mbox{as}\quad k\to\infty,
\end{equation}
therefore, by the first part of (\ref{b2}),
$$
\diam(\L_n\cup\G_n)\to0\quad\mbox{ as }\ k\to\infty\quad\mbox{ for }\ \L_n\cup\G_n\subset B(x_k,r_k).
$$
In particular, there is $k_2$ such that for $k>k_2$ and $\L_n\cup\G_n\subset B(x_k,r_k)$ we get
$$
\diam(\L_n\cup\G_n)<2^{-2N-1}\r,
$$
 and the set $\L_n\cup\G_n$  can intersect with only one of the  spherical shells
 $$
 B(a+\tau_k,2^{-j+1}\r)\setminus B(a+\tau_k,2^{-j}\r),\ j=1,\dots,2N+1.
 $$
On the other hand, by (\ref{b1}),
$$
\#\{n:\,(\L_n\cup\G_n)\subset B(a+\tau_k,2\r)\}\le 2N,
$$
 hence there is $m=m(k),\,1\le m\le 2N+1,$ such that
$$
(\L_n\cup\G_n)\cap[B(a+\tau_k,2^{-m+1}\r)\setminus B(a+\tau_k,2^{-m}\r)]=\emptyset\quad\mbox{ if }\ \L_n\cup\G_n\subset B(a+\tau_k,2\r).
$$
Consequently, the  sets $\L_n,\,\G_n$ are either both simultaneously subsets of $B(a+\tau_k,2^{-m}\r)$ and
$$
  \p\left(\frac{\l-a-\tau_k}{2^{-m}\r}\right)=\p\left(\frac{\g-a-\tau_k}{2^{-m}\r}\right)=1,\quad\forall\,\l\in\L_n,\ \forall\,\g\in\G_n,
$$
or are not subsets of $B(a+\tau_k,2^{-m+1}\r)$, and
$$
  \p\left(\frac{\l-a-\tau_k}{2^{-m}\r}\right)=\p\left(\frac{\g-a-\tau_k}{2^{-m}\r}\right)=0\quad\forall\,\l\in\L_n,\ \forall\,\g\in\G_n.
$$
Hence,
$$
 H_m(a+\tau_k)=\sum_{\l\in\L\cap B(a+\tau_k,2^{-m}\r)}\mu(\l)-\sum_{\g\in\G\cap B(a+\tau_k,2^{-m}\r)}\nu(\g)=\sum_{n:\L_n\subset B(a+\tau_k,2^{-m}\r)}[\mu(\L_n)-\nu(\G_n)].
  $$
By (\ref{n}) and the second part of (\ref{b2}), we obtain  $\mu(\L_n)-\nu(\G_n)\to0$ as $k\to\infty$. Since a number of members in the last sum is an most $N$, we get that the sequence $H_{m(k)}(a+\tau_k)$  tends to zero as $k\to\infty$, which contradicts to (\ref{e}). \bs
\medskip

{\bf Proof of Theorem \ref{T3}}.
 Assume the contrary $f\not\equiv g$. Then there is $a\in\supp f$ such that either $a\in\G$ and $p_{a,j^0}\neq q_{a,j^0}$ for some $j^0\in(\N\cup\{0\})^d$, or $a\not\in\G$ and $p_{a,j^0}\neq0$ for some $j^0\in(\N\cup\{0\})^d$.
In the latter case put $q_{a,j^0}=0$. Take $\p\in S(\R^d)$ such that $0\le\p\le1,\ \supp\p\subset B(0,2)$, and $\p(x)=1$ for $x\in B(0,1)$. Since $\L$ and $\G$ are discrete, it follows that there are no points of $\L\cup\G$ other than $a$ in the ball $B(a,2\r)$ for some $\r\le 1/2$. Put
$$
\psi(x)=\p(x/\r)x_1^{j_1^0}x_2^{j_2^0}\dots x_d^{j_d^0}/j_1^0!j_2^0!\dots j^0_d!,\qquad x=(x_1,\dots,x_d)\in\R^d.
$$
Clearly, $(D^{j^0}\psi)(0)=1$ and $(D^j\psi)(0)=0$ for $j\neq j^0$. Therefore,  $f(\psi_a)=p_{a,j^0}\neq g(\psi_a)$.
Using (\ref{f1}) and applying Lemma \ref{L2}, we obtain that the functions $f(\psi_t),\  g(\psi_t)$ are almost periodic, and the function $H(t)=f(\psi_t)-g(\psi_t)$ too.
Moreover, $H(a)=p_{a,j^0}- q_{a,j^0}\neq0$. Set $\e=|H(a)|/2$. Denote by $\hT$ the set of all $\e$-almost period of the functions $H(t)$. We get
\begin{equation}\label{g}
|H(a+\tau)|>\e\quad \forall\tau\in\hT.
\end{equation}
Since $\hT$ is relatively dense, it follows from (\ref{E}) that for every $k>k_1$ there is $\tau_k\in\hT$ such that $B(x_k,r_k)\supset B(a+\tau_k,\r)$.
Set
$$
   M_k=\{n:\,(\L_n\cup\G_n)\cap B(a+\tau_k,2\r)\neq\emptyset\}.
$$
By (\ref{b1}), $M_k$ is finite, hence,
\begin{equation}\label{f}
 \min M_k\to\infty \quad\mbox{ as }k\to\infty.
\end{equation}

Furthermore, for sufficiently large $k$
\begin{equation}\label{h}
f(\psi_{a+\tau_k})-g(\psi_{a+\tau_k})=\sum_{n\in M_k}\sum_j\left[\sum_{\l\in\L_n} p_{\l,j}(D^j\psi)(\l-a-\tau_k)-\sum_{\g\in\G_n}q_{\g,j}(D^j\psi)(\g-a-\tau_k)\right].
\end{equation}
All derivatives of $\psi$ are uniformly continuous, hence it follows from the first part of (\ref{f3}) and (\ref{f4}) that for some fixed points $b_n\in\L_n\cup\G_n$ and each $j$
$$
   \sum_{\l\in\L_n}p_{\l,j}\left[(D^j\psi)(\l-a-\tau_k)-(D^j\psi)(b_n-a-\tau_k)\right]\to0 \quad \mbox{ as }\ n\to\infty,
$$
and
$$
  \sum_{\g\in\G_n} q_{\g,j}\left[(D^j\psi)(\g-a-\tau_k)-(D^j\psi)(b_n-a-\tau_k)\right]\to0 \quad \mbox{ as }\ n\to\infty.
$$
 All derivatives of $\psi$ are uniformly bounded, hence it follows (\ref{f3}) that for each $j$
$$
   \left[\sum_{\l\in\L_n}p_{\l,j}-\sum_{\g\in\G_n} q_{\g,j}\right](D^j\psi)(b_n-a-\tau_k)\to0 \quad \mbox{ as }\ n\to\infty.
  $$
  Hence,
$$
   \sum_{\l\in\L_n}p_{\l,j}(D^j\psi)(\l-a-\tau_k)-\sum_{\g\in\G_n} q_{\g,j}(D^j\psi)(\g-a-\tau_k)\to0 \quad \mbox{ as }\ n\to\infty.
$$
 Note that the number of values of $j$ does not exceed $(m+1)^d$ and $\# M_k\le 2N$.   Thus we obtain from (\ref{f}) and (\ref{h})
 $$
    |H(a+\tau_k)|=|f(\psi_{a+\tau_k})-g(\psi_{a+\tau_k})|\to0 \quad \mbox{ as }\ k\to\infty,
 $$
  which contradicts to (\ref{g}). \bs
\medskip

Theorem \ref{T4} follows from Theorems \ref{T2} and \ref{T3}, if only we change the Fourier transform to the inverse Fourier transform.  Theorems \ref{T1} and \ref{T1a} follow from Theorem \ref{T2},
and  Theorem \ref{T5} follows from Theorem \ref{T4}.
\medskip

{\bf Question}. It can be proved that Theorem \ref{T2} is valid for the case when $\mu,\nu$ are sums of atomic and absolutely continuous components. If it is valid for arbitrary almost periodic
 measures?
\medskip

 The author is grateful to Mikhail Sodin for drawing the attention of the author to the article \cite{KS}, and to  Pavel Kurasov for his interest in the research of the author and useful discussion.

\end{document}